\documentclass[preprint,12pt]{elsarticle} 
\usepackage[english]{babel}

\usepackage{enumerate}          
\usepackage{amsmath}
\usepackage{amssymb}    
\usepackage{graphicx}     
\usepackage{epsfig}           
\usepackage{enumerate}

\journal{xxx}  
\begin{document} 
\begin{abstract}
In this paper, we are interested in the relation between the solutions of the control system $\dot x=f(x,u)$  and the solutions of its (potentially unknown) perturbation $\dot x=f(x,u)+w(x,t).$ Under the assumption that the linear part of the unperturbed system at the point $(0,0)$ is controllable and that disturbance $w(x,t)$ is sufficiently small, there exists a state-feedback controller of the form $u=-Kx$ such that the perturbed system preserves the local asymptotic stability of the zero solution of unperturbed system. The main result of this paper  gives the sufficient conditions, more specifically, the relations between the important parameters of the system, to ensure this property and at the same time provides the method for calculating the lower bound of region of attraction. Moreover, we obtain a nontrivial extension of the classical result of H.~K.~Khalil regarding asymptotic behavior of the (uncontrolled) perturbed systems whose nominal part is exponentially asymptotically stable at the origin $x=0.$
\end{abstract}


\title{On the local asymptotic stabilization of the nonlinear systems with small time-varying perturbations by state-feedback control}
\author{R.~Vrabel}
\ead{robert.vrabel@stuba.sk}
\address{Slovak University of Technology in Bratislava, Institute of Applied Informatics, Automation and Mechatronics,  Bottova 25,  917 01 Trnava, Slovakia}

\newtheorem{thm}{Theorem}
\newtheorem{lem}[thm]{Lemma}
\newtheorem{defi}[thm]{Definition}
\newdefinition{rmk}{Remark}
\newdefinition{ex}{Example}
\newproof{pf}{Proof}
\newproof{pot1}{Proof of Theorem \ref{thm1}}

\pagestyle{headings}

\maketitle

\section[Introduction]{Introduction}

In this paper we will study the nonlinear control system
\begin{equation}\label{system:unperturbed}
\dot x=f(x,u), \ t\geq t_0,
\end{equation}
with the state-feedback of the form $u=-Kx$ and for which is assumed that $x=0$ is its solution, that is $f(0,0)=0.$ It is well-known, that if the pair $(A,B),$ where $A=f_x(0,0)$  and $B=f_u(0,0)$ are the corresponding Jacobian matrices of the vector field $f(x,u)$ with respect to the state and input variables, respectively, and evaluated at $(0,0)$  is controllable, then the LTI system $\dot x=(A-BK)x$ is in some neighborhood of the origin topologically equivalent (and preserving the parametrization by time) to the system $\dot x=f(x,-Kx),$ provided that the eigenvalues of the matrix $A-BK$ do not lie on the imaginary axis. The precise statement about this property gives the Hartman-Grobman theorem (see, e.~g. \cite[p.~120]{Perko}).  

Now, if the system (\ref{system:unperturbed}) is subjected to certain small, time-varying perturbation, the above property concerning the solutions near the origin may or may not remain true. A more precise formulation
of this problem is as follows: if $x=0$ is locally asymptotically stable (l.~a.~s.) in the sense of Lyapunov solution for unperturbed problem (\ref{system:unperturbed}) with $u=-Kx,$ and if $w(x,t)$ is sufficiently small, give conditions on $f,$ $w$ and the state-feedback gain matrix $K$  so that origin is  l.~a.~s. for the perturbed system
\begin{equation}\label{system:perturbed}
\dot x=f(x,u)+w(x,t), \ u=-Kx, \ t\geq t_0.
\end{equation}
Note that $x=0$ may not be its solution. However, it is possible to overcome this problem by utilizing the following definition of asymptotic stability of the origin which is a natural generalization of its usual form  used in the Lyapunov stability theory.
\begin{defi}[compare with \cite{Strauss}, Definitions~2.2-2.4] \label{definition_elas}
We say, that the origin $x=0$ is {\bf eventually l.~a.~s.} for the perturbed closed-loop system (\ref{system:perturbed}) if for every $\tilde\varepsilon>0$ there exists $t^*(\tilde\varepsilon)\geq0$ and $\tilde\delta(\tilde\varepsilon)>0$ such that if $|x(t_0)|<\tilde\delta$ then $|x(t)|<\tilde\varepsilon$ for all $t\geq t_0+t^*,$ and, moreover, $|x(t)|\rightarrow 0^+$ with $t\rightarrow \infty.$ A stronger form of the eventual local asymptotic stability is if $|x(t)|\leq \alpha_1|x(t_0)|e^{\alpha_2(t-t_0)}$ for all $t\geq t_0+t^*,$ where $\alpha_1>0,$ $\alpha_2<0$ and $|x(t_0)|<\tilde\delta,$ that is, the solution of (\ref{system:perturbed}) converges for $t\rightarrow \infty$ to the origin exponentially.
\end{defi}
The origin is l.~a.~s. in the usual sense, see, e.~g. \cite[Definition~4.1, p.~112]{Khalil}, if one can choose $t^*=0.$ This implies at the same time that $x=0$ is a solution of (\ref{system:perturbed}), as is proved in \cite[Lemma~2.7]{Strauss}.

Due to the highly complicated nature of the nonlinear control systems, finding an analytical relation between the governing equations (represented by a vector field $f$), the gain matrix $K$ and perturbation $w$ ensuring (at least) local asymptotic stability is not easy. 

For comparison, for the linear systems $\dot x=A(t)x$ with a continuous on $[t_0,\infty)$ matrix function $A(t),$ the sufficient conditions for eventual global asymptotic stability of perturbed system (in Definition~\ref{definition_elas}, $|x(t)|\rightarrow 0^+$ for any initial state $x(t_0)\in\mathbb{R}^n$) are 
\begin{itemize}
\item[(a)] $\|\Phi(t,\tau)\|\leq\alpha_A e^{\beta_A(t-\tau)},$ for some $\alpha_A>0,$ $\beta_A<0$ and for all $t,$ $\tau$ such that $t\geq\tau,$ where $\Phi$  is a state transition matrix of $\dot x=A(t)x$ and
\item[(b)] for all $x\in\mathbb{R}^n$ and $t\geq t_0$ is $|w(x,t)|\leq g_w(t)\rightarrow 0^+$ as $t\rightarrow\infty$ 
\end{itemize}
as follows from the inequality holding for all $t\geq t_0,$
\[
|x(t)|\leq \|\Phi(t,t_0)\||x(t_0)|+\int_{t_0}^t \|\Phi(t,\tau)\||w(x(\tau),\tau)|d\tau
\]
\[
\leq \alpha_A e^{\beta_A(t-t_0)}|x(t_0)|+\alpha_A\frac{\int_{t_0}^t e^{-\beta_A\tau}g_w(\tau)d\tau}{e^{-\beta_At}}
\]
and L'Hospital's Rule applied on the second term on the right--hand side of the inequality above.
Obviously, the condition (a) is satisfied if $A$ is a constant matrix and every eigenvalue of $A$ has a strictly negative real part.

In contrast,  an example is known (\cite{Strauss}, Theorem~C and Example~8.2) in which $f(x,t)$ is uniformly continuous and locally Lipschitz on $\mathbb{R}\times[0,\infty)$ and all solutions of unperturbed system approach zero exponentially and monotonically as $t\rightarrow \infty$ but origin is not attracting for perturbed system $\dot x=f(x,t)+\gamma_w(x)e^{-\beta_w t}$ for each $\beta_w\geq0$ and each continuous function $\gamma_w$ satisfying $\gamma_w(x)>0$ for $x>0.$ In the light of these findings, every result regarding stability of the perturbed nonlinear systems has its value and significance.
 
There are two useful methods for studying the qualitative behavior of solutions of the nonlinear systems: {\bf{(i)}} the Lyapunov's second method and its various extensions (such as the La Salle's invariance principle),  see, e.~g. \cite{Amato}, \cite{Chesi1}, \cite{Chesi2}, \cite{ChiangHirschWu}, \cite{Najafi}, \cite{Ran}, \cite{Topcu}, \cite{Valmorbida}, \cite{Vannelli} and {\bf{(ii)}} the use of variation of constants formula where the behavior of solutions  of a perturbed system is determined in terms of the behavior of solutions of an unperturbed system (\cite{Gonzales}).  This idea is essential in the proof of Lyapunov's indirect method in Corollary~2.43 in \cite[p.~160]{Chicone}. Neither  of these mentioned articles, and as far as we know, nor any other, has dealt with the systematic study of the analytical calculation method for region of attraction for perturbed systems if the perturbing term is time-varying. 

In principle, there are two possible approaches to the problem under consideration. One approach is to set conditions on $f$ and $K,$  and find out what kind of admissible perturbations $w$ preserve local asymptotic stability. The second approach is reversed, to set the growth conditions on the perturbation $w$ that will be allowed, and find out which control systems (by manipulating their parameters) will have their asymptotic stability preserved by all such $w.$  

Following the first approach,  by combining the controllability theory, the method of variation constants and three-functions variant of the Gronwall-Bellman inequality (Lemma~\ref{Gronwall}), we will establish the sufficient conditions on the perturbation $w$ (in the form of growth constraints) to be preserved the local asymptotic stability of the zero solution of the state-feedback control system $\dot x=f(x,-Kx).$ 

\section{Notations and Assumptions}

Let $x, f, w$ and $u$  are  $n-$ and  $m-$dimensional column vectors, respectively.  We shall always assume that $f$ is continuous in $(x,u)$ and twice continuously differentiable with respect to the components of $x$ and $u$ for  $|x|<\infty,$ $|u|<\infty,$ that $f(0,0) = 0,$  that $w$ is continuous in $(x, t)$ for $|x|<\infty,$ $t_0\leq t<\infty,$  and that the domain of existence of trajectories for the control systems under consideration is the interval $[t_0,\infty)$ for every initial state $x(t_0)\in\mathbb{R}^n.$

The $n\times n$ matrix $f_x(0,0)$ and  $n\times m$ matrix $f_u(0,0)$ are the Jacobian matrices of the vector function $f(x,u)$ with respect to the variables $x$ and $u,$ respectively, and evaluated at $(x,u)=(0,0),$  $K$ is an $m\times n$ constant gain matrix and an upper dot indicates a time derivative. We denote by $|\cdot|$ the Euclidean norm and by $\|\cdot\|$ a matrix norm induced by the Euclidean norm of vectors, $\|A\|=\max_{|x|=1}|Ax|.$   The real part of a complex number $z$ is denoted by $\Re(z)$ and  the superscript $'\,{\mathsf{T}}\,'$ is used to indicate transpose operator.

In the following section, with the aim of weakening conservatism in estimating, we formulate an important technical result as a consequence of the slight variant of the Gronwall-Bellman inequality given in \cite[p.~56]{Bellman} which allows us more subtle to estimate the qualitative behavior of the solutions of perturbed system in the proof of the main result of the present paper, Theorem~\ref{thm1}.
 
\section[main result]{Gronwall-Bellman type inequality for three functions and the main result}

\begin{lem}\label{Gronwall}
If $U(t),V(t),W(t)\geq0$ are the continuous functions for all $t\geq t_0,$ if $C$ is a positive constant, and if
\begin{equation*}
U(t)\leq C+\int\limits_{t_0}^t\left(U(\tau)V(\tau)+W(\tau)\right)d\tau
\end{equation*}
then
\begin{equation*}
U(t)\leq e^{\int_{t_0}^t V(\tau)d\tau}\left[C+\int_{t_0}^t W(\tau)d\tau\right]
\end{equation*}
for all $t\geq t_0.$ 
\end{lem}
\begin{pf} 
The proof of this lemma follows by the fact that $C+\int_{t_0}^t W(\tau)d\tau$ as a function of the variable $t$ is positive, monotonic and nondecreasing, together with Lemma~1 given in \cite{Pachpatte}.
\end{pf}

Now we formulate our result.

\begin{thm}\label{thm1}
Let us consider the control system (\ref{system:perturbed}), namely,
\[
\dot x=f(x,u)+w(x,t), \ u=-Kx, \ t\geq t_0.
\] 
Assume that
\begin{itemize}
\item[1.] $f(0,0)=0;$
\item[2.] the LTI pair $(A, B),$ where $A=f_x(0,0),$ $B=f_u(0,0),$ is controllable; 
\item[3.] for all $t\geq t_0$ is $|w(x,t)|\leq ce^{\gamma(t-t_0)},$ $c\geq0,$ $\gamma<0.$
\end{itemize}
Let $\lambda_{m}\triangleq\max\left\{\Re(\lambda_i),\, i=1,\dots,n \right\},$ $\lambda_i$'s are the eigenvalues of the matrix $A_{cl}=A-BK,$ $\lambda_i\neq\lambda_j$ and let $\lambda_{m}<0.$

Let $\eta\triangleq\|T\|\|T^{-1}\|,$ where $T$ is a similarity transformation for which $\Lambda=T^{-1}A_{cl}T=diag(\lambda_i).$ 

Let $\Gamma_0>0$ and $\varepsilon_0>0$ are such that $|R_1(x,-Kx)|\leq\Gamma_0|x|$ for $|x|\leq\varepsilon_0;$ here $R_1(x,u)$ means the Taylor expansion remainder of the function $f(x,u)$ expanded at the point $(x,u)=(0,0).$

Further assume that for the desired eigenvalues $\lambda_i,$ $i=1,\dots,n$ of the matrix $A_{cl}=A-BK$ there exists $\Gamma_0>0$ such that
\begin{itemize}
\item[($\ast$)] $\lambda_{m}<-\eta\Gamma_0(1+\|K\|).$
\end{itemize}
Then 
\begin{itemize}
\item[a)] the origin $x=0$ is eventually l.~a.~s. for (\ref{system:perturbed}) in the sense of Definition~\ref{definition_elas};
\item[b)] for all $x(t_0)$ satisfying $|x(t_0)|\leq\delta$ is $|x(t)|\leq\varepsilon_0$ for all $t\geq t_0,$ and moreover, $|x(t)|\rightarrow 0^+$ (exponentially) for $t\rightarrow\infty,$ where $\delta$ is a positive solution of the equation
\begin{equation}\label{eq:delta}
\eta\left(\delta+{\frac{ c}{\lambda_{m}-\gamma}}\right)=\varepsilon_0,
\end{equation} 
that always exists for sufficiently small values of $c\geq0,$ and $\gamma,$ $\gamma<\lambda_{m}<0;$ 
\item[c)] The radius of the lower bound of region of attraction $\delta=\delta(\varepsilon_0, \lambda_m, c, \gamma)\leq\varepsilon_0/\eta$ with equality only if $c=0,$ that is, for $w(x,t)\equiv 0.$
\end{itemize}
\end{thm}
Before proving Theorem~\ref{thm1} we demonstrate its applicability in approximating the region of attraction of control system.
 
\begin{ex}\label{example1}
Consider the system
\begin{equation}\label{example_system}
\dot x_1=x_1^3+x_2^2u_1+x_2,\ \dot x_2=u_1+ w_2(x_1,x_2,t),\ t\geq0.
\end{equation}
We first verify that the linearization of this system is controllable,
\[
A=f_x(0,0)=\left(\begin{array}{cc}
 0 \ & 1   \\
 0 \ & 0
\end{array} \right)\ \mathrm{and}\ B=f_u(0,0)=\left(\begin{array}{c}
 0    \\
 1
\end{array} \right)
\]
and thus the controllability matrix $(B \ AB)$ is regular. Therefore the linear part of original system is controllable.

Now, let us calculate the relation between $\Gamma_0$ and $\varepsilon_0.$ The Taylor remainder $R_1(x,u)$ consists of the second-order partial derivatives of the function $f$ and in our case we have the estimate
\[
|R_1(x,u)|\leq\frac12\Big[6|x_1\|x|^2+2(2|x_2|)|x_1\|u_1|+2|u_1\|x|^2 \Big].
\]
Since $|x_1\|x_2|\leq|x|^2/2$ and $|u_1|\leq\|K\||x|,$
\[
|R_1(x,-Kx)|\leq\Big[ \max\{3,\|K\|\} +\|K\| \Big]|x|^3.
\]
Thus, the inequality $|R_1(x,-Kx)|\leq\Gamma_0|x|$ in Theorem~\ref{thm1} holds if 
\[
|x|\leq\sqrt{\Gamma_0/(\max\{3,\|K\|\} +\|K\|)}\ (=\varepsilon_0).
\]
Now, let $\lambda_1=-0.5$ and $\lambda_2=-0.75$ are the desired closed-loop eigenvalues. This implies that the gain matrix $K=(0.3750\ 1.2500)$ and the associated norms are $\|K\|=1.3050$ and $\eta=11.0902.$ So, from the Assumption~($\ast$) of Theorem~\ref{thm1}, $\Gamma_0<-\lambda_m/\eta(1+\|K\|)=0.0196$ and the corresponding limiting value $\varepsilon_0=0.0674.$

If we choose $c=0.001$ and $\gamma=-10$ we obtain from the equation (\ref{eq:delta}) that $\delta=0.0059.$  Theorem~\ref{thm1} guarantees that the solution of the problem under consideration with the state-feedback $u_1=-K(x_1,x_2)^{\mathsf{T}}$ starting with $|x(0)|<0.0059$ converges to zero with $t\rightarrow \infty,$ and $|x(t)|<0.0674$ for all $t\geq0.$
\begin{figure}
   \centerline{
    \hbox{
     \psfig{file=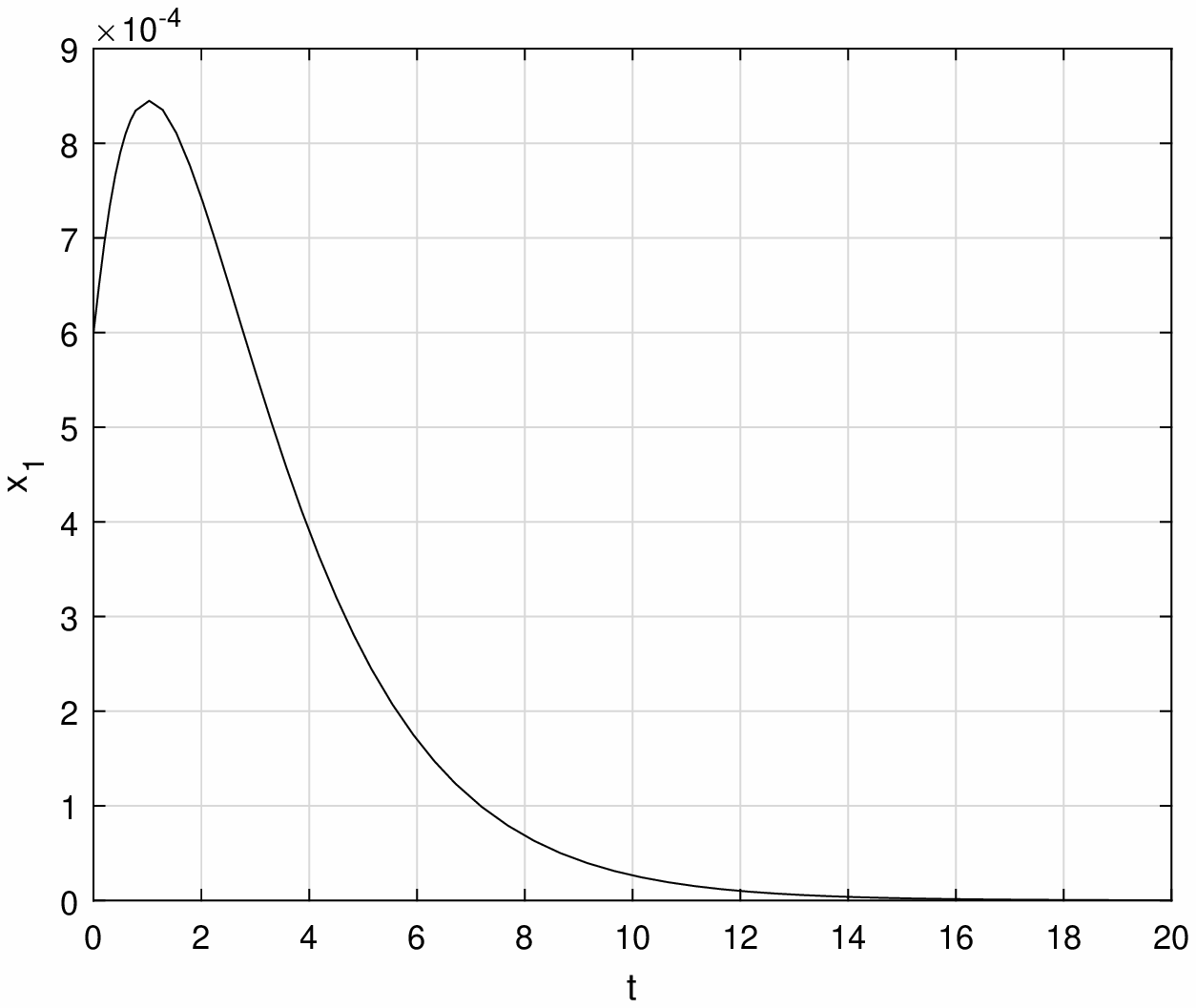,width=6cm, clip=}
     \hspace{1.cm}
     \psfig{file=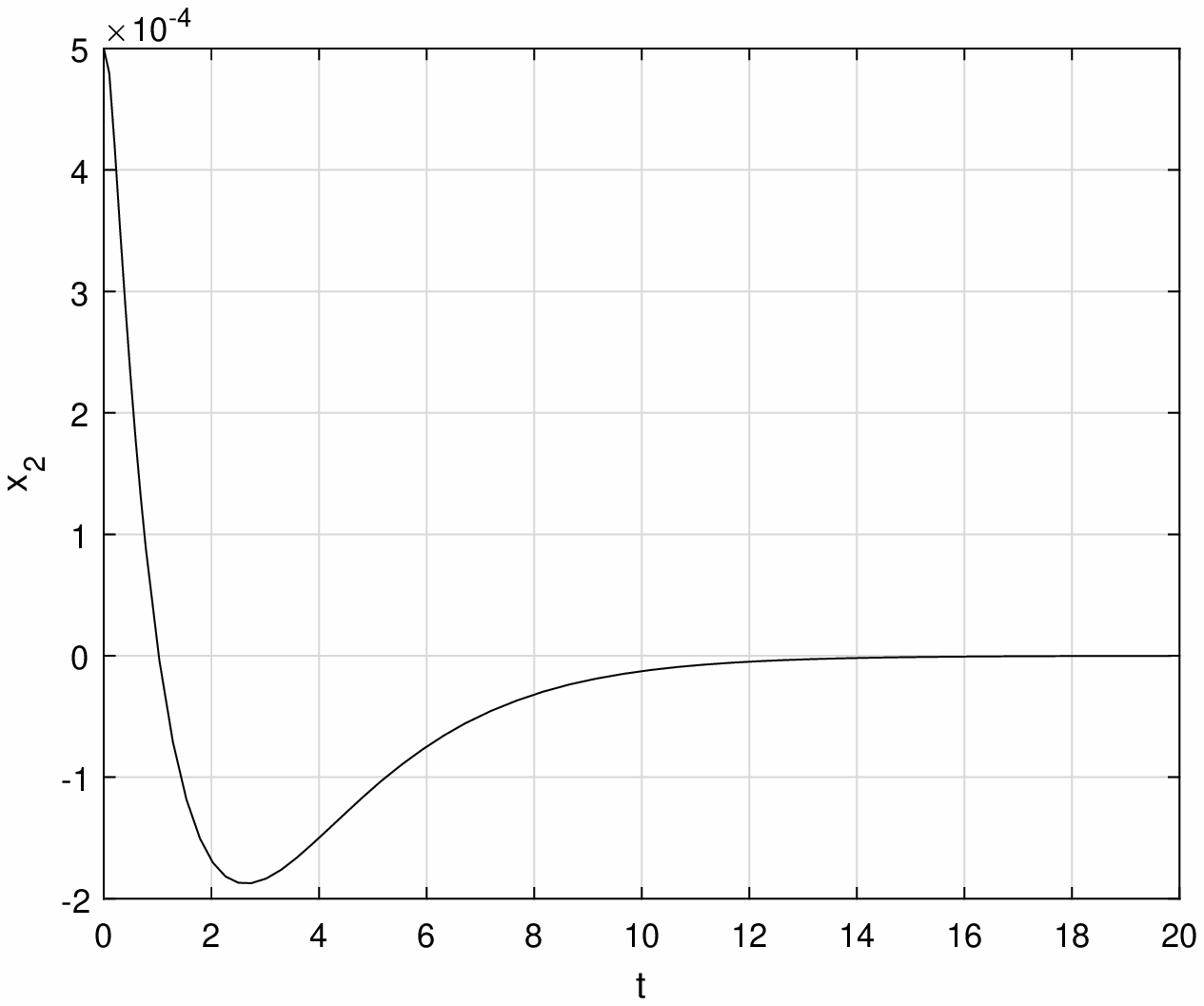,width=6cm,clip=}
    }
   }
\caption{The solution of control system $\dot x_1=x_1^3+x_2^2u_1+x_2,$ $\dot x_2=u_1+ w_2(x_1,x_2,t)$ with the state-feedback $u_1=-Kx,$ $K=(0.3750 \ 1.2500)$ for perturbation $w_2(x_1,x_2,t)=0.001e^{-10t}\cos(x_1+x_2)$ and with the initial state $x(0)=(6,5)^{\mathsf{T}}\times10^{-4}.$ }
\label{fig_example}
\end{figure} 
Figure~\ref{fig_example} demonstrates the effectiveness of the proposed controller locally asymptotically stabilizing the perturbed system in the neighborhood of the origin.
\end{ex}
\begin{pot1} 
The system $\dot x=f(x,u)+w(x,t)$  may be rewritten, after the Taylor expansion of $f$ at $(x,u)=(0,0),$ in the form
\begin{equation*}
\dot x = Ax+Bu+\tilde R(x,u,t),\ A=f_x(0,0),\ B=f_u(0,0)
\end{equation*}
and $\tilde R(x,u,t)=R_1(x,u)+w(x,t),$ where $R_1$ is the Taylor remainder in Lagrange's form. Since $R_1(x,u)$ is clearly $O(|(x,u)|^2)$ as $|(x,u)|\rightarrow 0^+,$ we can find the constants $\varepsilon_0>0$ and  $\Gamma_0>0$ such that $|R_1(x,-Kx)|\leq \Gamma_0\left(|x|+|-Kx|\right)$ for $|x|\leq\varepsilon_0.$  From the last inequality and the Assumption~3 of Theorem~\ref{thm1} we have the estimate
\begin{equation}\label{estimate:tildeR}
|\tilde R(x,-Kx,t)|\leq \Gamma_0(1+\|K\|)|x|+ce^{\gamma(t-t_0)}\ \mathrm{for}\ |x|<\varepsilon_0,\ t\geq t_0.
\end{equation}
Identifying $\tilde R$ as an inhomogeneous term,  every solution of (\ref{system:perturbed}) satisfies a nonlinear integral equation
\begin{equation*}
x(t)=e^{A_{cl}(t-t_0)}x(t_0)+\int_{t_0}^t e^{A_{cl}(t-\tau)}\tilde R(x(\tau),-Kx(\tau),\tau)d\tau.
\end{equation*}
From the properties of the matrix exponentials and the used operator norm we have that
\[
\|e^{A_{cl}(t-t_0)}\|=\|e^{T\Lambda T^{-1}(t-t_0)}\|=\|Te^{\Lambda(t-t_0)}T^{-1}\|
\]
\[
\leq\underbrace{\|T\|\|T^{-1}\|}_{\triangleq\eta}\|e^{\Lambda(t-t_0)}\|=\eta e^{\lambda_{m}(t-t_0)},
\]
where $\eta$ is no less than one \cite[p.~101]{Rugh}, and so
\[
|x(t)|\leq\eta e^{\lambda_{m}(t-t_0)}|x(t_0)|+\int_{t_0}^t \eta e^{\lambda_{m}(t-\tau)}|\tilde R(x(\tau),-Kx(\tau),\tau)|d\tau.
\]
Then, multiplying by $e^{-\lambda_{m}t}$ and using (\ref{estimate:tildeR}), we get
\[
e^{-\lambda_{m}t}|x(t)|\leq\eta e^{-\lambda_{m}t_0}|x(t_0)|+\int_{t_0}^t \eta e^{-\lambda_{m}\tau}\Big[\underbrace{\Gamma_0(1+\|K\|)}_{\triangleq\Theta }|x(\tau)|+ce^{\gamma(\tau-t_0)}\Big]d\tau.
\]
Now, if we associate  $e^{-\lambda_{m}t}|x(t)|$ with $U(t),$ $\eta\Theta$ with $V(t),$ $\eta ce^{-\lambda_{m}t}e^{\gamma(t-t_0)}$ with $W(t),$ $\eta e^{-\lambda_{m}t_0}|x(t_0)|$ with $C,$ and using Lemma~\ref{Gronwall} 
\[
e^{-\lambda_{m}t}|x(t)|\leq e^{\int_{t_0}^t \eta\Theta d\tau}\Bigg[\eta e^{-\lambda_{m}t_0}|x(t_0)|+\eta c{\int_{t_0}^t e^{-\lambda_{m}\tau}e^{\gamma(\tau-t_0)}d\tau}\Bigg],
\]
and after an algebraic manipulation
\[
|x(t)|\leq\eta e^{(\eta\Theta+\lambda_m)(t-t_0)}\Bigg[|x(t_0)|+{ce^{(\lambda_m-\gamma) t_0}\int_{t_0}^t e^{(\gamma-\lambda_m)\tau}d\tau}\Bigg].
\]
For $\gamma-\lambda_m<0,$ the integral in the square bracket is convergent for $t\rightarrow \infty$ and
\[
\int_{t_0}^t e^{(\gamma-\lambda_m)\tau}d\tau\leq\int_{t_0}^{\infty} e^{(\gamma-\lambda_m)\tau}d\tau\leq\frac{e^{(\gamma-\lambda_m)t_0}}{\lambda_m-\gamma}<\infty.
\]  
Using this, we get the final estimate for the solutions $x(t)$ of perturbed problem,
\begin{equation}\label{final_estimate}
|x(t)|\leq \eta e^{(\eta\Theta+\lambda_m)(t-t_0)}\Bigg[ |x(t_0)|+\left({\frac{c}{\lambda_m-\gamma}}\right)\Bigg],\ t\geq t_0,
\end{equation}
proving that the origin $x=0$ is eventually l.~a.~s. in the sense of Definition~\ref{definition_elas}
if $\lambda_{m}<-\eta\Theta=-\eta\Gamma_0(1+\|K\|)$ and $\gamma<\lambda_m.$ The time shift $t^*>0$ in the mentioned definition compensates the second term in the square brackets of the above estimate; $t^*=0$ if and only if $c=0.$ Moreover, because the expression in the square brackets is constant, the solution $x(t)$ of perturbed problem converges exponentially fast to zero and $|x(t)|\leq\varepsilon_0$  for all initial states satisfying $|x(t_0)|\leq\delta,$ where $\delta$ is a positive solution of (\ref{eq:delta}) which always exists if $c$ is sufficiently small.  This completes the proof of Theorem~\ref{thm1}. 
\end{pot1}
\begin{rmk}
The important additional value of the just proved theorem is in providing method of calculating the lower bound of  attractivity region for system, illustrated in Example~\ref{example1}.
\end{rmk}

\begin{rmk}
The statement of Theorem~\ref{thm1}  one may generalize all the above steps in the proof to the control systems with perturbations of the form 
\[
|w(x,t)|\leq\sigma|x| +ce^{\gamma(t-t_0)},\ \sigma\geq0,\ c\geq0,\ \gamma<0,\ t\geq t_0,
\] 
where the condition $(*)$ is replaced by the following one:
\[
\lambda_{m}<-\eta\bigg[\Gamma_0(1+\|K\|)+\sigma\bigg].
\]
Thus we obtain a nontrivial generalization of the classical result of Khalil \cite[Lemma~9.4, p.~352]{Khalil} regarding asymptotic behavior of the (uncontrolled) perturbed systems whose nominal part is exponentially asymptotically stable at origin $x=0$ giving, by using the converse Lyapunov theorem, the boundedness of solutions of perturbed systems only.
\end{rmk}

\section*{Conclusions}
This paper dealt with the local asymptotic stability of the origin of perturbed nonlinear control systems  $\dot x =f(x,u)+w(x,t),$ $u=-Kx,$ $f(0,0)=0,$ $t\geq t_0,$ where $K$ is the state-feedback gain matrix and the perturbing term $w(x,t)$ $\bigg[\left\vert w(x,t)\right\vert\leq\sigma|x|+ ce^{\gamma(t-t_0)},$ $\sigma\geq0,$ $c\geq0,$ $\gamma<0,$ $x\in\mathbb{R}^n,$ $t\geq t_0\bigg]$ represents the (potentially unknown) time-varying disturbance and/or unmodeled dynamics of system. Because $x=0$ may not be solution of perturbed system, a generalization of the notion "local asymptotic stability" was introduced. It is shown, that by appropriate choice of the growth constraints imposed on the disturbance term $w$ through the parameters $c,$ $\gamma$ in the dependence on the gain matrix $K,$ the perturbed system preserves the local asymptotic stability property of the unperturbed system $\dot x =f(x,u),$ $u=-Kx$ in the neighborhood of its equilibrium $x=0.$ For this purpose, in Theorem~\ref{thm1} we defined the relations between the important parameters of the system and moreover, as additional value, together with the proof of theorem we derived the method for calculating the lower bound of the region of attraction which was applied in Example~\ref{example1}. 

As one of the main used analytical tools, we have used three-functions variant of Gronwall-Bellman inequality, which produced a closed formula for the closed-loop solutions of perturbed system on the interval $[t_0,\infty),$ namely, that
\[
|x(t)|\leq \tilde\alpha_1 e^{\alpha_2(t-t_0)}\Bigg[ |x(t_0)|+{\alpha_3}\Bigg],
\]
where the constants $\tilde\alpha_1\geq1,$ $\alpha_2<0,$ $\alpha_3\geq0$ [$\tilde\alpha_1=\tilde\alpha_1(K),$ $\alpha_2=\alpha_2(f,K,\sigma)$  and $\alpha_3=\alpha_3(c,\gamma,K),$ $\alpha_3=0$ iff $c=0$]. This result extends the now classical result of Khalil about asymptotics of perturbed systems whose nominal (unperturbed) part is exponentially stable.

\end{document}